\documentclass[a4paper,12pt,leqno]{article}

\usepackage{graphicx}
\usepackage{oldgerm}
\usepackage{theorem}

\usepackage{amsfonts}
\usepackage{latexsym}
\usepackage{epsfig}
\usepackage{amssymb}
\usepackage{graphicx}
\usepackage{oldgerm}
\usepackage{theorem}

\def\0{{\bf 0}}

\def\g{g_{N,T}}

\def\x{{\bf x}}
\def\y{{\bf y}}

\def\mg{\mathfrak{g}_{N,T}}

\def\R{\mathbb{R}}

\def\mX{\mathfrak{X}}
\def\X{{\bf X}}

\def\tS{\widetilde{S}}
\def\tI{\widetilde{I}}

\def\RV{\mathbb{R}^{N}_{<}}

\def\cN{{\cal N}_N}

\def\Pf{{\rm Pf}}
\def\vtheta{\mib{\theta}}

\usepackage{theorem}

\newcommand{\mib}[1]{\mbox{\boldmath $#1$}}

\theorembodyfont{\itshape}
\newtheorem{thm}{Theorem}

\begin{document}
\pagestyle{plain}
\vskip 0.3cm

\begin{center}
{\bf \Large{Non-Colliding System of Brownian Particles}}
\vskip 0.3cm

{\bf \Large{as Pfaffian Process}}

\vskip 0.5cm

Makoto Katori
\footnote{Electronic mail: katori@phys.chuo-u.ac.jp}\\
{\it Department of Physics,
Faculty of Science and Engineering,\\
Chuo University, 
Kasuga, Bunkyo-ku, 
Tokyo 112-8551, Japan}

\end{center}

\vskip 0.5cm

In the paper \cite{KNT04} we studied the 
temporally inhomogeneous system of
non-colliding Brownian motions and proved that 
multi-time correlation functions are generally given
by the quaternion determinants in the sense of Dyson and Mehta.
In this report we give another proof of the equivalent
statement using Fredholm determinant and
Fredholm pfaffian, and claim that the present system
is a typical example of pfaffian processes.

\vskip 0.5cm

\setcounter{equation}{0}
\section{Non-Colliding Brownian Motions}

By virtue of the Karlin-McGregor formula 
\cite{KM59a, KM59b}, the transition density 
of the absorbing Brownian motion in a Weyl chamber
$$
\RV=\Big\{\x=(x_{1}, x_{2}, \dots, x_{N}) \in \R^{N};
x_{1} < x_{2} < \cdots < x_{N} \Big\},
$$ 
is given by
$$
f_{N}(t; \x, \y)=\det_{1 \leq i, j \leq N} \Big[
p_t (x_{i}, y_{j}) \Big], \quad
\x, \y \in \RV, \quad t \in [0, \infty),
$$
where $p_{t}$ is the heat-kernel given by
$$
p_t (x,y)= \frac{1}{\sqrt{2 \pi t}} e^{-(x-y)^2/2t}.
$$
The integral
$$
\cN(t; \x)=\int_{\RV} f_{N}(t; \x, \y) d \y,
$$
where $d\y=\prod_{i=1}^{N} dy_{i}$,
gives the probability that a Brownian motion 
starting from $\x\in \RV$ does not hit the boundary of $\RV$
up to time $t>0$.

For a given $T > 0$, we define
\begin{equation}
\g(s, \x; t, \y)
=\frac{f_{N}(t-s; \x, \y) \cN(T-t;\y)}
{\cN(T-s; \x)}
\label{eqn:gNT}
\end{equation}
for $0 \leq s \leq t \leq T, \x, \y \in \RV$.
It can be regarded as the transition probability density from the
state $\x \in \RV $ at time $s$ 
to the state $\y \in \RV$ at time $t$.
In \cite{KT02, KT03a} it was shown that
as $|\x|\to 0$, $g_{N}^{T}(0,\x; t,\y)$ converges to 
\begin{equation}
\g(0, {\bf 0}; t, \y)
\equiv C(N,T,t)  h_{N}(\y) \prod_{i=1}^N p_t(0,y_i) \cN(T-t, \y),
\label{eqn:gNT0}
\end{equation}
where
$C(N,T,t)=\pi^{N/2}\left(\prod_{j=1}^{N} \Gamma(j/2)\right)^{-1}
T^{N(N-1)/4} t^{-N(N-1)/2}$, and
$$
  h_{N}(\x) = \prod_{1 \leq i<j \leq N} (x_{j}-x_{i})
  = \det_{1 \leq i, j \leq N} \Big[ x_{j}^{i-1} \Big].
$$
The $N$ particle system of non-colliding Brownian motions
$\X(t)$
all started from the origin at time $0$,
{\it i.e.} $\X(0)=\0=(0,0, \cdots, 0)$,
and conditioned not 
to collide with each other in a time interval $(0,T]$ is defined by the
process associated with the transition probability density $\g$
given by (\ref{eqn:gNT}) and (\ref{eqn:gNT0}).
This process is temporally inhomogeneous
and it was obtained as a diffusion scaling limit of
the vicious walker model in \cite{KT03a}.
Figure 1 illustrates the process $\X(t)$.

\begin{figure}[htbp]
\includegraphics[width=.6\linewidth]{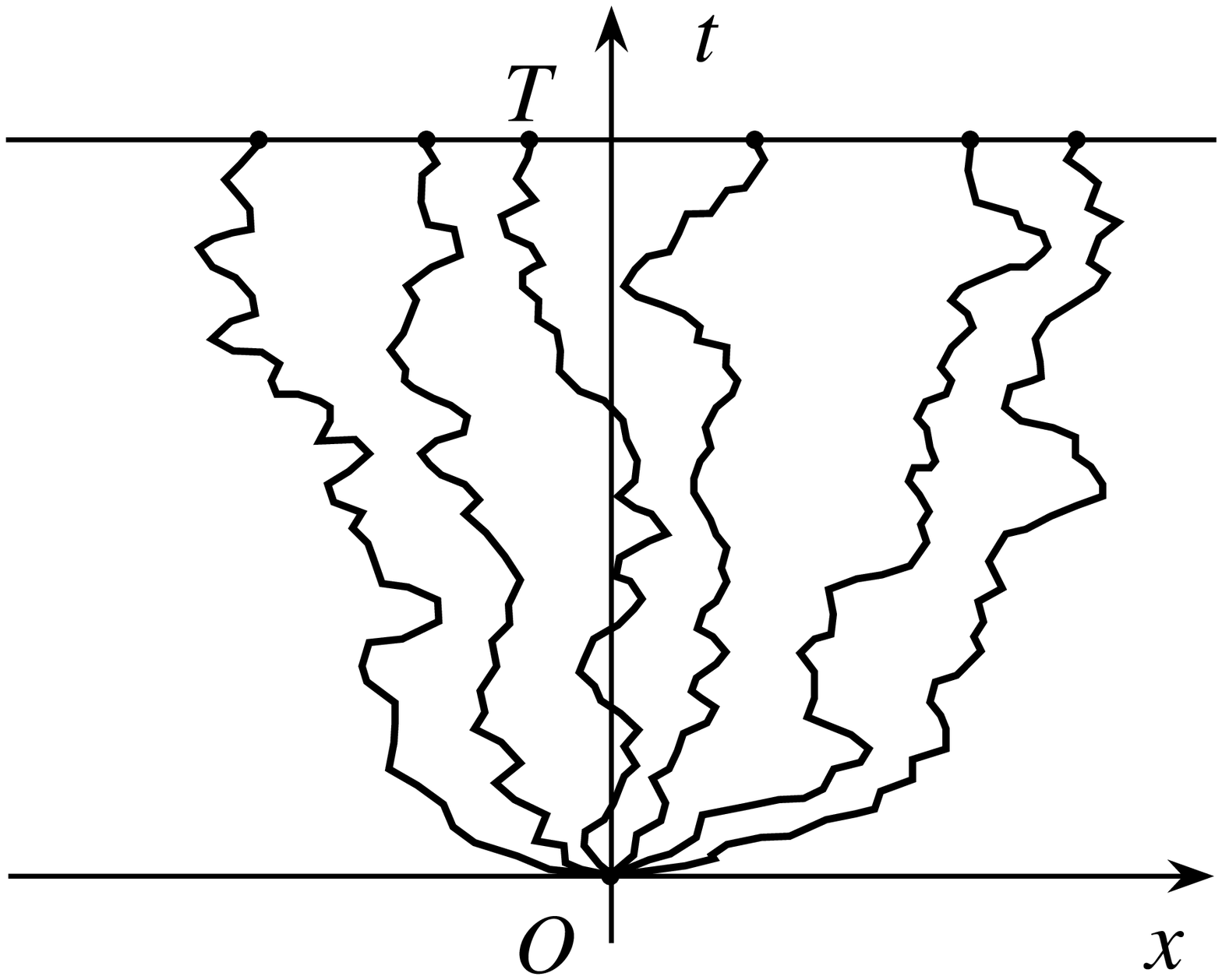}
\caption{Process $\X(t), t \in [0, T]$, with $\X(0)=\0$.
\label{fig:star}}
\end{figure}

Let $\mX$ be the space of countable subset $\xi$ of $\R$
satisfying $\sharp (\xi \cap K) < \infty$
for any compact subset $K$.
For $\x = (x_1, x_2, \dots, x_n)\in \bigcup_{\ell=1}^{\infty} \R^\ell$,
we denote $\{x_i\}_{i=1}^n \in \mX$ simply by $\{\x\}$.
The diffusion process $\{\X(t)\}, t \in [0,T]$ on the set $\mX$
is defined by the transition probability density
$\mg(s, \{\x\} ; t, \{\y\})$, $0 \leq s \leq t \leq T$:
$$
\mg(s, \{\x\} ; t, \{\y \})=
\left\{
\begin{array}{ll}
\g(s, \x; t, \y),
& \mbox{if} \ s>0,  \ \sharp \{\x\} = \sharp \{\y\} =N,
\\
\g(0, {\bf 0}; t, \y),
& \mbox{if} \ s=0, \ \{\x\} =\{ 0 \}, \ \sharp \{\y\}=N,
\\
0, 
& \mbox{otherwise},
\end{array}\right.
$$
where $\x$ and $\y$ are the elements of $\RV$. 
For the given time interval $[0,T]$, we consider the $M$
intermediate times
$t_{0}=0 < t_{1} < t_{2} < \cdots < t_{M} < t_{M+1}=T$.
The multi-time transition density
of the process $\{\X(t)\}$ is given by
$$
\mg(t_{0}, \{\x^{(0)}\} ;t_{1}, \{\x^{(1)}\};
\dots; t_{M+1}, \{\x^{(M+1)}\}) 
=\prod_{\mu=0}^{M} \mg(t_{\mu}, \{\x^{(\mu)}\}; 
t_{\mu+1}, \{\x^{(\mu+1)} \}).
$$
Here we set
$t_{0}=0$, $t_{M+1}=T$ and $\{\x^{(0)} \} \equiv \{ 0\}$.
From (\ref{eqn:gNT}) and (\ref{eqn:gNT0})
we have 
\begin{eqnarray}
\label{eqn:mg}
&& \mg(0, \{ 0 \} ;t_{1}, \{\x^{(1)} \};
\dots; t_{M+1}, \{\x^{(M+1)} \})
\\
&& \quad
= C(N,T,t_{1}) h_{N}\left(\x^{(1)} \right) {\rm sgn}
\left( h_{N} \left(\x^{(M+1)} \right) \right)
\nonumber\\
&& \qquad \times
\prod_{i=1}^N p_{t_1}\left(0, x_i^{(1)}\right)
\prod_{\mu=1}^{M} 
\det_{1 \leq i, j \leq N} \left[
p_{t_{\mu+1}-t_{\mu}} \left(x_{i}^{(\mu)}, x_{j}^{(\mu+1)} 
\right) \right].
\nonumber
\end{eqnarray}
For $\x^{(\mu)} \in \RV$, $1\leq \mu \leq M+1$,
and $N'=1,2,\dots, N$, we write
$\x^{(\mu)}_{N'} = \left(x_1^{(\mu)}, x_2^{(\mu)}, \dots, 
x_{N'}^{(\mu)}\right)$.
For a sequence $\{N_\mu \}_{\mu=1}^{M+1}$ of positive integers 
less than or equal to $N$,
we define the multi-time 
correlation function by
\begin{eqnarray}
\label{eqn:corr}
&&
\rho_{N,T} \left(t_{1}, \{\x^{(1)}_{N_1}\}; t_2, \{\x^{(2)}_{N_2}\}; 
\dots; t_{M+1}, \{\x^{(M+1)}_{N_{M+1}}\} \right) \\
&&\qquad =
\int \limits_{\prod_{\mu=1}^{M+1} \R^{N-N_{\mu}}}
\prod_{\mu=1}^{M+1}
\frac{1}{(N-N_{\mu})!}\prod_{i=N_{\mu}+1}^{N} dx_{i}^{(\mu)} \nonumber\\
&& \hskip 4cm
\mg \left(0, \{0 \}; t_{1}, \{\x^{(1)}_{N}\}; t_2, \{\x^{(2)}_{N}\};
\dots; t_{M+1}, \{\x^{(M+1)}_{N}\} \right).
\nonumber
\end{eqnarray}

Let $C_{0}(\R)$ be the set of all continuous real functions
with compact supports. 
For ${\bf f}=(f_{1}, f_{2}, \cdots, f_{M+1}) \in C_{0}(\R)^{M+1}$,
and $\vtheta=(\theta_{1}, \theta_{2}, \cdots, \theta_{M+1})
\in \R^{M+1}$,
the multi-time characteristic function is defined
for the process $\X(t)=(X_{1}(t), X_{2}(t), \cdots, X_{N}(t))$ as
$$
\Psi_{N,T} ( {\bf f}; \vtheta)
= {\bf E}_{N,T} \left[
\exp \left\{ \sqrt{-1} \sum_{\mu=1}^{M+1}
\theta_{\mu} \sum_{i_{\mu}=1}^{N} 
f_{\mu}(X_{i_{\mu}}(t_{\mu})) \right\} \right],
$$
where ${\bf E}_{N,T}[\, \cdot \, ]$ denotes the expectation determined
by $\mg$. 
Let 
$$
\chi_{\mu}(x)=e^{\sqrt{-1} \theta_{\mu} f_{\mu}(x)}-1, 
\quad 1 \leq \mu \leq M+1,
$$
then by the definition of
multi-time correlation function (\ref{eqn:corr}), we have
\begin{eqnarray}
\Psi_{N,T} ( {\bf f}; \vtheta )
&=& \sum_{N_{1}=0}^{N} \sum_{N_{2}=0}^{N} \cdots
\sum_{N_{M+1}=0}^{N}
\int_{\R^{N_{1}}} d \x_{N_{1}}^{(1)}
\int_{\R^{N_{2}}} d \x_{N_{2}}^{(2)} \cdots
\int_{\R^{N_{M+1}}} d \x_{N_{M+1}}^{(M+1)} \nonumber\\
\label{eqn:corr2}
&& \prod_{\mu=1}^{M+1} \prod_{i_{\mu}=1}^{N_{\mu}} 
\chi_{\mu}(x_{i_{\mu}}^{(\mu)})
\rho_{N,T} \left(t_{1}, \{\x_{N_{1}}^{(1)} \}; 
t_{2}, \{\x_{N_{2}}^{(2)} \};
\cdots ; t_{M+1}, \{\x_{N_{M+1}}^{(M+1)}\} \right).
\end{eqnarray}

\setcounter{equation}{0}
\section{Fredholm Pfaffian Representation
of Characteristic Function}

Let
\begin{eqnarray}
Z_{N,T}[\chi] &=& \left(\frac{1}{N !}\right)^{M+1}
\int_{\R^{N(M+1)}} \prod_{\mu=1}^{M+1} d \x^{(\mu)} \,
\det_{1 \leq i, j \leq N} \left[
M_{i-1}(x_{j}^{(1)}) p_{t_{1}}(0, x_{j}^{(1)})
(1+ \chi_{1}(x_{j}^{(1)})) \right] \nonumber\\
&& \times
\prod_{\mu=1}^{M} \det_{1 \leq i, j \leq N}
\left[ p_{t_{\mu+1}-t_{\mu}} (x_{i}^{(\mu)}, x_{j}^{(\mu+1)})
(1+\chi_{\mu+1}(x_{j}^{(\mu+1)})) \right]
{\rm sgn}\left(h_{N}(\x_{N}^{(M+1)}) \right).
\nonumber
\end{eqnarray}
Here $M_{i}(x)$ is an arbitrary polynomial
of $x$ with degree $i$ in the form
$M_{i}(x)=b_{i} x^{i}+ \cdots$ 
with a constant $b_{i}$
for $i=0,1,2, \cdots$, and thus
$h_{N}(\x)=\det_{1 \leq i, j \leq N} \Big[
M_{i-1}(x_{j}) \Big]/ \prod_{k=1}^{N} b_{k-1}$. 
Then (\ref{eqn:mg}) gives
\begin{equation}
\label{eqn:Psi1}
\Psi_{N,T} ( {\bf f}; \vtheta )=
\frac{Z_{N,T}[\chi]}{Z_{N,T}[0]},
\end{equation}
where
$Z_{N,T}[0]$ is obtained from $Z_{N,T}[\chi]$ by setting
$\chi_{\mu}(x)=0$ for all $\mu$.

By the well-known formula
$$
\int_{\RV} d \x \,
\det_{1 \leq i, j \leq N} \Bigg[\phi_{i}(x_{j})\Bigg]
\det_{1 \leq i, j \leq N} \Bigg[ \bar{\phi}_{i}(x_{j}) \Bigg]
= \det_{1 \leq i, j \leq N} \left[
\int_{\R} dx \, \phi_{i}(x) \bar{\phi}_{j}(x) \right]
$$
for square integrable continuous functions
$\phi_{i}, \bar{\phi}_{i}, 1 \leq i \leq N$,
we have
\begin{eqnarray}
Z_{N,T}[\chi] &=& \int_{\RV} d \y \, 
\det_{1 \leq i, j \leq N} \left[
\int_{\R^{M+1}} \prod_{\mu=1}^{M+1} dx^{(\mu)} \,
\Bigg\{ M_{i-1}(x^{(1)}) p_{t_{1}}(0, x^{(1)})
(1+ \chi_{1}(x^{(1)})) \Bigg\} \right. \nonumber\\
&& \quad \left. \times \prod_{\mu=1}^{M} \Bigg\{
p_{t_{\mu+1}-t_{\mu}}(x^{(\mu)}, x^{(\mu+1)}) 
(1+\chi_{\mu+1}(x^{(\mu+1)}))
\Bigg\}
p_{T-t_{M+1}}(x^{(M+1)}, y_{j}) \right],
\nonumber
\end{eqnarray}
where 
$p_{T-t_{M+1}}(x,y)=p_{0}(x,y)=\delta(x-y)$.

For simplicity of expressions, we will assume that 
the number of particles $N$ is even from now on. We use the
formula \cite{Bruijn55}
$$
\int_{\RV} d\y \,
\det_{1 \leq i, j \leq N} \Bigg[
\phi_{i}(y_{j}) \Bigg]
= \Pf_{1 \leq i, j \leq N} \left[
\int_{\R} dy \int_{\R} d \tilde{y} \,
{\rm sgn}(\tilde{y}-y) \phi_{i}(y) \phi_{j}(\tilde{y})
\right]
$$
for integrable continuous functions
$\phi_{i}, 1 \leq i \leq N$.
Here,
for an integer $n$ and an antisymmetric $2n \times 2n$ matrix
$A=(a_{ij})$, the pfaffian is defined as
$$
\Pf(A) = \Pf_{1 \leq i < j \leq 2n} \Big[ \ a_{ij} \ \Big] 
=\frac{1}{n !} \sum_{\sigma}
{\rm sgn}(\sigma) a_{\sigma(1) \sigma(2)} a_{\sigma(3) \sigma(4)} 
\cdots a_{\sigma(2n-1) \sigma(2n)},
$$
where the summation is extended over all permutations $\sigma$
of $(1,2,\dots, 2n)$ with restriction
$\sigma(2k-1) < \sigma(2k), k=1,2,\dots, n$.
Since $(\Pf(A))^2=\det A$ for any antisymmetric
$2n \times 2n$ matrix $A$, we have
\begin{eqnarray}
\Bigg( Z_{N,T}[\chi] \Bigg)^2
&=& \det_{1 \leq i, j \leq N} \left[
\int dy \int d \tilde{y} \, {\rm sgn}(\tilde{y}-y) \right. 
\nonumber\\
&& \ \times 
\int_{\R^{M+1}} \prod_{\mu=1}^{M+1} dx^{(\mu)}
\Bigg\{ M_{i-1}(x^{(1)}) p_{t_{1}}(0, x^{(1)})
(1+ \chi_{1}(x^{(1)})) \Bigg\} \nonumber\\
&& \quad \times \prod_{\mu=1}^{M} \Bigg\{
p_{t_{\mu+1}-t_{\mu}}(x^{(\mu)}, x^{(\mu+1)}) 
(1+\chi_{\mu+1}(x^{(\mu+1)}))
\Bigg\}
p_{T-t_{M+1}}(x^{(M+1)}, y) \nonumber\\
&& \ \times 
\int_{\R^{M+1}} \prod_{\mu=1}^{M+1} d\tilde{x}^{(\mu)}
\Bigg\{ M_{j-1}(\tilde{x}^{(1)}) 
p_{t_{1}}(0, \tilde{x}^{(1)})
(1+ \chi_{1}(\tilde{x}^{(1)})) \Bigg\} \nonumber\\
&& \quad \left. \times \prod_{\mu=1}^{M} \Bigg\{
p_{t_{\mu+1}-t_{\mu}}(\tilde{x}^{(\mu)}, \tilde{x}^{(\mu+1)}) 
(1+\chi_{\mu+1}(\tilde{x}^{(\mu+1)}))
\Bigg\}
p_{T-t_{M+1}}(\tilde{x}^{(M+1)}, \tilde{y}) \right] \nonumber\\
&=& \det_{1 \leq i, j \leq N} \Bigg[
(A_{0})_{ij}+(A_{1})_{ij}+(A_{2})_{ij}+(A_{3})_{ij} \Bigg]
\nonumber
\end{eqnarray}
with
\begin{eqnarray}
(A_{0})_{ij} &=& \int dy \int d \tilde{y} \, {\rm sgn}(\tilde{y}-y) \,
\int dx M_{i-1}(x) p_{t_{1}}(0,x) p_{T-t_{1}}(x,y) \nonumber\\
&& \qquad \times \int d \tilde{x} M_{j-1}(\tilde{x})
p_{t_{1}}(0,\tilde{x}) p_{T-t_{1}}(\tilde{x}, \tilde{y}),
\nonumber\\
(A_{1})_{ij} &=& \sum_{\ell=1}^{M+1} 
\sum_{1 \leq \mu_{1} < \mu_{2} < \cdots < \mu_{\ell} \leq M+1}
\int dy \int d \tilde{y} \, {\rm sgn}(\tilde{y}-y) \nonumber\\
&& \quad \times \int_{\R^{\ell}} \prod_{k=1}^{\ell} d x^{(\mu_{k})}
\int d x \, M_{i-1}(x) p_{t_{1}}(0, x) 
p_{t_{\mu_{1}}-t_{1}}(x, x^{(\mu_{1})}) \nonumber\\
&& \qquad \times \prod_{k=1}^{\ell-1} \Bigg\{
\chi_{k}(x^{(\mu_{k})}) p_{t_{\mu_{k+1}}-t_{\mu_{k}}}
(x^{(\mu_{k})}, x^{(\mu_{k+1})})
\Bigg\} \chi_{\ell}(x^{(\mu_{\ell})}) p_{T-t_{\mu_{\ell}}}
(x^{(\mu_{\ell})}, y)
\nonumber\\
&& \quad \times \int d \tilde{x} M_{j-1}(\tilde{x}) 
p_{t_{1}}(0, \tilde{x}) p_{T-t_{1}}(\tilde{x}, \tilde{y}),
\nonumber\\
(A_{2})_{ij} &=& \sum_{\ell=1}^{M+1} 
\sum_{1 \leq \mu_{1} < \mu_{2} < \cdots < \mu_{\ell} \leq M+1}
\int dy \int d \tilde{y} \, {\rm sgn}(\tilde{y}-y) \nonumber\\
&& \quad \times \int d x M_{i-1}(x) 
p_{t_{1}}(0, x) p_{T-t_{1}}(x, y),
\nonumber\\
&& \quad \times \int_{\R^{\ell}} \prod_{k=1}^{\ell} 
d \tilde{x}^{(\mu_{k})}
\int d \tilde{x} \, M_{j-1}(\tilde{x}) 
p_{t_{1}}(0, \tilde{x}) 
p_{t_{\mu_{1}}-t_{1}}(\tilde{x}, \tilde{x}^{(\mu_{1})}) \nonumber\\
&& \qquad \times \prod_{k=1}^{\ell-1} \Bigg\{
\chi_{k}(\tilde{x}^{(\mu_{k})}) 
p_{t_{\mu_{k+1}}-t_{\mu_{k}}}(\tilde{x}^{(\mu_{k})}, \tilde{x}
^{(\mu_{k+1})})
\Bigg\} \chi_{\ell}(\tilde{x}^{(\mu_{\ell})})
p_{T-t_{\mu_{\ell}}}(\tilde{x}^{(\mu_{\ell})}, \tilde{y}),
\nonumber\\
(A_{3})_{ij} &=& \sum_{\ell=1}^{M+1} \sum_{m=1}^{M+1} 
\sum_{1 \leq \mu_{1} < \mu_{2} < \cdots < \mu_{\ell} \leq M+1}
\sum_{1 \leq \nu_{1} < \nu_{2} < \cdots < \nu_{m} \leq M+1}
\int dy \int d \tilde{y} \, {\rm sgn}(\tilde{y}-y) \nonumber\\
&& \quad \times \int_{\R^{\ell}} \prod_{k=1}^{\ell} d x^{(\mu_{k})}
\int d x \, M_{i-1}(x) p_{t_{1}}(0, x) 
p_{t_{\mu_{1}}-t_{1}}(x, x^{(\mu_{1})}) \nonumber\\
&& \qquad \times \prod_{k=1}^{\ell-1} \Bigg\{
\chi_{k}(x^{(\mu_{k})}) p_{t_{\mu_{k+1}}-t_{\mu_{k}}}
(x^{(\mu_{k})}, x^{(\mu_{k+1})})
\Bigg\} \chi_{\ell}(x^{(\mu_{\ell})}) p_{T-t_{\mu_{\ell}}}
(x^{(\mu_{\ell})}, y)
\nonumber\\
&& \quad \times \int_{\R^{\ell}} \prod_{n=1}^{m} 
d \tilde{x}^{(\nu_{m})}
\int d \tilde{x} \, M_{j-1}(\tilde{x}) 
p_{t_{1}}(0, \tilde{x}) 
p_{t_{\nu_{1}}-t_{1}}(\tilde{x}, \tilde{x}^{(\nu_{1})}) \nonumber\\
&& \qquad \times \prod_{n=1}^{m-1} \Bigg\{
\chi_{n}(\tilde{x}^{(\nu_{n})}) 
p_{t_{\nu_{n+1}}-t_{\nu_{n}}}(\tilde{x}^{(\nu_{n})}, 
\tilde{x}^{(\nu_{n+1})})
\Bigg\} \chi_{m}(\tilde{x}^{(\nu_{m})})
p_{T-t_{\nu_{m}}}(\tilde{x}^{(\nu_{m})}, \tilde{y}),
\nonumber
\end{eqnarray}
where we have used the Chapman-Kolmogorov equation
for the heat-kernel
$$
\int  dy \, p_{t-s}(x,y) p_{u-t}(y,z)=p_{u-s}(x,z),
\quad 0 < s < t < u, \, x, y \in \R,
$$
and the fact that $p_{0}(x,y)=\lim_{t \to 0} p_{t}(x,y)
=\delta(x-y)$.

We consider a vector space ${\cal V}$ with the orthonormal
basis $\Big\{ |\mu,x \rangle \Big\}_{
\mu=1,2, \cdots, M+1, x \in \R}$ which satisfy
$$
 \langle \mu, x | \nu, y \rangle =\delta_{\mu \nu} \delta(x-y),
$$
$\mu,\nu =1,2, \cdots, M+1, x, y \in \R$.
We introduce the operators $\hat{J}, \hat{p},
\hat{p}_{+}, \hat{p}_{-}$ and $\hat{\chi}$ 
acting on ${\cal V}$ as follows
\begin{eqnarray}
&& \langle \mu, x | \hat{J} | \nu, y \rangle
= {\bf 1}_{(\mu=\nu=M+1)} {\rm sgn}(y-x), \nonumber\\
&& \langle \mu, x | \hat{p} | \nu, y \rangle
= p_{|t_{\nu}-t_{\mu}|}(x,y), \nonumber\\
\label{eqn:p+-}
&& \langle \mu, x | \hat{p}_{+} | \nu, y \rangle
= p_{t_{\nu}-t_{\mu}}(x,y) {\bf 1}_{(\mu<\nu)}
= \langle \nu, y| \hat{p}_{-} | \mu, x \rangle, \\
&& \langle \mu, x | \hat{\chi} | \nu, y \rangle
= \chi_{\mu}(x) \delta_{\mu \nu} \delta(x-y),
\nonumber
\end{eqnarray}
where
${\bf 1}_{(\omega)}$ is the indicator function:
${\bf 1}_{(\omega)}=1$ if $\omega$ is satisfied
and ${\bf 1}_{(\omega)}=0$ otherwise,
and we will use the convention
$$
\langle \mu, x | \hat{A} | \nu, y \rangle
\langle \nu, y | \hat{B} | \rho, z)
=\sum_{\nu=1}^{M+1} \int_{\R} dy \,
A(\mu,x;\nu,y)B(\nu,y;\rho,z)
=\langle \mu,x | \hat{A}\hat{B}|\rho,z \rangle
$$
for operators $\hat{A}, \hat{B}$ with
$\langle \mu, x | \hat{A} | \nu, y \rangle=A(\mu,x;\nu,y),
\langle \mu, x | \hat{B} | \nu, y \rangle=B(\mu,x;\nu,y)$.

Consider another basis $\Big\{ |i \rangle \, ; \,
i=1,2, \cdots \Big\}$ in ${\cal V}$ and we assume that
transformation matrix between 
the two bases is given by
\begin{equation}
\label{eqn:R}
\langle i | \mu,x \rangle 
= \langle \mu, x | i \rangle
=\int dy M_{i-1}(y) p_{t_{1}}(0,y)
p_{t_{\mu}-t_{1}}(y,x),
\end{equation}
$i=1,2, \cdots, \mu=1,2, \cdots, M+1, x \in \R$.
Then the quantity
\begin{eqnarray}
\label{eqn:psi}
&& \int dy \, M_{i-1}(y) p_{t_{1}}(0,y) p_{t_{\mu}-t_{1}}(y,x)
\nonumber\\
&&+\sum_{m=1}^{\mu-1} \sum_{1 \leq \mu_{1} < \mu_{2} < 
\cdots < \mu_{m} < \mu}
\int_{\R^{m}} \prod_{j=1}^{m} dx^{(\mu_{j})} \,
\int dy \, M_{i-1}(y) p_{t_{1}}(0,y) 
p_{t_{\mu_{1}}-t_{1}}(y,x^{(\mu_{1})})
\nonumber\\
&& \qquad \quad \times \prod_{k=1}^{m-1} \Bigg\{
\chi_{\mu_{k}}(x^{(\mu_{k})}) p_{t_{\mu_{k+1}}-t_{\mu_{k}}}
(x^{(\mu_{k})}, x^{(\mu_{k+1})}) \Bigg\}
\chi_{\mu_{m}}(x^{(\mu_{m})}) 
p_{t_{\mu}-t_{\mu_{m}}}(x^{(\mu_{m})},x),
\nonumber
\end{eqnarray}
for $\mu=1,2, \cdots, M+1, x \in \R, i=1,2, \cdots $,
can be written as
\begin{eqnarray}
&& \langle i | \mu, x \rangle
+ \sum_{m \geq 1} \langle i | \mu_{1}, x^{(\mu_{1})} \rangle
\langle \mu_{1}, x^{(\mu_{1})} | \hat{\chi} \hat{p}_{+} | 
\mu_{2}, x^{(\mu_{2})} \rangle
\nonumber\\
&& \qquad \qquad \qquad
\cdots 
\langle \mu_{m-1}, x^{(\mu_{m-1})} | \hat{\chi} \hat{p}_{+} | 
\mu_{m}, x^{(\mu_{m})} \rangle
\langle \mu_{m}, x^{(\mu_{m})} | \hat{\chi} \hat{p}_{+} | 
\mu, x \rangle
\nonumber\\
&=& \langle i | \mu,x \rangle
+ \sum_{m \geq 1} \langle i| 
( \hat{\chi} \hat{p}_{+} )^{m} | \mu, x \rangle
\nonumber\\
&=& \langle i | \frac{1}{1- \hat{\chi} \hat{p}_{+}}
| \mu, x \rangle. \nonumber
\end{eqnarray}
It is also expressed as \,
$\displaystyle{
   \langle \mu, x | \frac{1}{1-\hat{p}_{-} \hat{\chi}} | i \rangle.
}$

It should be noted that the basis
$\Big\{ |i \rangle \, ; \,
i=1,2, \cdots \Big\}$ 
is in general not orthonormal.
Here we introduce an operator $\hat{\delta}$ such that
\begin{equation}
\label{eqn:delta1}
\langle i | \hat{\delta} | j \rangle
= \langle j | \hat{\delta} | i \rangle
= \delta_{ij}, \quad
i,j=1,2, \cdots.
\end{equation}
We will use the following convention
$$
\langle i | \hat{A} | j \rangle
\langle j | \hat{B} | \mu, x \rangle
= \sum_{j=1}^{\infty} A_{ij} B_{j}^{(\mu)}(x)
$$
for
$
A_{ij} = \langle i | \hat{A} | j \rangle
\equiv \langle i | \mu, x \rangle
\langle \mu, x | \hat{A} | \nu, y \rangle
\langle \nu, y | j \rangle, 
B_{j}^{(\mu)}(x) =
\langle j | \hat{B} | \mu, x \rangle
\equiv \langle j | \nu, y \rangle
\langle \nu, y | \hat{B} | \mu, x \rangle, 
$
but we will not write it as
$\langle i | \hat{A} \hat{B} | \mu, x \rangle$,
since 
$\Big\{ |i \rangle \, ; \,
i=1,2, \cdots \Big\}$
is in general not a complete basis.
By this basis, any operator $\hat{A}$ on ${\cal V}$
may have a semi-infinite matrix representation
$A=\Big(\langle i | \hat{A} | j \rangle 
\Big)_{i,j =1,2, \cdots}$.
If the matrix $A$ representing an operator $\hat{A}$
is invertible, we define the operator
$\hat{A}^{\bigtriangleup}$ such that its matrix representation
is the inverse of $A$;
\begin{equation}
\label{eqn:Delta1}
  \Big( \langle i | \hat{A}^{\bigtriangleup} | j \rangle
  \Big)_{i,j=1,2, \cdots} = A^{-1}.
\end{equation}
In other words,
$$
\langle i | \hat{A} |j \rangle
\langle j | \hat{A}^{\bigtriangleup} | k \rangle
=\langle i | \hat{A}^{\bigtriangleup} |j \rangle
\langle j | \hat{A} | k \rangle
= \langle i | \hat{\delta} | k \rangle,
\qquad i,k=1,2, \cdots.
$$
For any given operator $\hat{A}$, the equality
$$
\langle i | \hat{A} | j \rangle
\langle j | \hat{A}^{\bigtriangleup} | k \rangle
\langle k | \hat{B} | \ell \rangle
= \langle i | \hat{B} | \ell \rangle.
$$
holds for arbitrary $i, \ell=1,2,\cdots$ and
$\hat{B}$. Then the equality
\begin{equation}
\label{eqn:equalityA}
\hat{A} |j \rangle \langle j |
\hat{A}^{\bigtriangleup} | k \rangle
\langle k | = 1
\end{equation}
should be established for each $\hat{A}$.
We will use this equality later.

Let ${\cal P}_{N}$ be a projection operator from the space
spanned by $\Big\{ |i \rangle \, ; \,
i=1,2, \cdots \Big\}$ to its $N$-dimensional subspace
spanned by
$\Big\{ |i \rangle \, ; \,
i=1,2, \cdots, N \Big\}$, and thus 
$$
\langle i | {\cal P}_{N} | \mu, x \rangle
= \langle \mu, x | {\cal P}_{N} | i \rangle
= \left\{
\begin{array}{ll}
\langle i |\mu, x \rangle,
& \mbox{if} \ 1 \leq i \leq N,
\\
& \\
0, 
& \mbox{otherwise}.
\end{array}\right.
$$
Then we have the following expressions 
for $(A_{\alpha})_{ij},
\alpha=0,1,2,3, i,j =1,2, \cdots, N$,
\begin{eqnarray}
\label{eqn:A0}
(A_{0})_{ij} &=& \langle i | {\cal P}_{N}| \mu, x \rangle
\langle \mu, x| \hat{J} | \nu, y \rangle
\langle \nu, y |{\cal P}_{N}| j \rangle \nonumber\\
&=& \langle i | {\cal P}_{N} \hat{J} {\cal P}_{N} |
j \rangle, \nonumber\\
(A_{1})_{ij} &=& \langle i |{\cal P}_{N}
\frac{1}{1-\hat{\chi}\hat{p}_{+}}\hat{\chi} | \mu, x \rangle
\langle \mu, x | \hat{p} \hat{J} {\cal P}_{N}
| j \rangle \nonumber\\
&=& \langle i |{\cal P}_{N} \hat{\chi}
\frac{1}{1-\hat{p}_{+} \hat{\chi}} | \mu, x \rangle
\langle \mu, x | \hat{p} \hat{J} {\cal P}_{N}
| j \rangle \nonumber\\
&=& \langle i | {\cal P}_{N}  \hat{\chi}
\frac{1}{1-\hat{p}_{+} \hat{\chi}}
\hat{p} \hat{J} {\cal P}_{N}| j \rangle, \nonumber\\
(A_{2})_{ij} &=& \langle i | {\cal P}_{N} \hat{J} \hat{p} 
| \mu, x \rangle
\langle \mu, x | \hat{\chi} 
\frac{1}{1-\hat{p}_{-}\hat{\chi}}{\cal P}_{N}
|j \rangle, \nonumber\\
&=& \langle i | {\cal P}_{N} \hat{J} \hat{p}
\hat{\chi} \frac{1}{1-\hat{p}_{-}\hat{\chi}}{\cal P}_{N} |
j \rangle, \nonumber\\
(A_{3})_{ij} &=& \langle i | {\cal P}_{N}
\frac{1}{1-\hat{\chi}\hat{p}_{+}}\hat{\chi}
| \mu, x \rangle 
\langle \mu, x | \hat{p}\hat{J} \hat{p} |\nu, y \rangle
\langle \nu, y | \hat{\chi} \frac{1}{1-\hat{p}_{-}\hat{\chi}} 
{\cal P}_{N}| j \rangle
\nonumber\\
&=& \langle i | {\cal P}_{N}\hat{\chi}
\frac{1}{1-\hat{p}_{+}\hat{\chi}}
| \mu, x \rangle 
\langle \mu, x | \hat{p}\hat{J} \hat{p} |\nu, y \rangle
\langle \nu, y | \hat{\chi} \frac{1}{1-\hat{p}_{-}\hat{\chi}} 
{\cal P}_{N}| j \rangle
\nonumber\\
&=& \langle i | {\cal P}_{N}\hat{\chi}
\frac{1}{1-\hat{p}_{+}\hat{\chi}}
 \hat{p}\hat{J} \hat{p} \hat{\chi}
 \frac{1}{1-\hat{p}_{-}\hat{\chi}} {\cal P}_{N} | j \rangle.
 \nonumber
\end{eqnarray}
That is, the matrices
$A_{\alpha}=\Big( (A_{\alpha})_{ij} \Big)_{i,j=1,2, \cdots, N},
\alpha=0,1,2,3,$ can be regarded as the
matrix representations in the basis
$\Big\{ |i \rangle \, ; \, i=1,2, \cdots \Big\}$
of the operators.

Since $\Big( Z_{N,T}[0] \Big)^2
=\det_{1 \leq i, j \leq N} \Big[ (A_{0})_{ij} \Big]$, 
(\ref{eqn:Psi1}) gives
\begin{equation}
\label{eqn:det1}
\Bigg\{ \Psi_{N,T} ( {\bf f}; \vtheta ) 
\Bigg\}^2
= \det_{1 \leq i, j \leq N} \Bigg[
\delta_{ij}+(A_{0}^{-1}A_{1})_{ij}+(A_{0}^{-1}A_{2})_{ij}
+(A_{0}^{-1}A_{3})_{ij} \Bigg].
\end{equation}
By our notation (\ref{eqn:Delta1}),
$$
(A_{0}^{-1})_{ij}=
\langle i | {\cal P}_{N} 
({\cal P}_{N} \hat{J} {\cal P}_{N})^{\bigtriangleup}
{\cal P}_{N} | j \rangle,
$$
since it satisfies the relation
\begin{eqnarray}
&& \langle i | {\cal P}_{N}
({\cal P}_{N} \hat{J}{\cal P}_{N})^{\bigtriangleup} 
{\cal P}_{N} | j \rangle
\langle j | {\cal P}_{N} \hat{J}{\cal P}_{N} | k \rangle
\nonumber\\
&=& {\bf 1}_{(1 \leq i, k \leq N)} 
\sum_{j=1}^{N} (A_{0}^{-1})_{ij} (A_{0})_{jk}
= {\bf 1}_{(1 \leq i, k \leq N)} \delta_{ik} \nonumber\\
&=& \langle i | {\cal P}_{N} \hat{\delta}
{\cal P}_{N} | j \rangle.
\nonumber
\end{eqnarray}
We will use the abbreviation
$$
   \hat{A}_{N}={\cal P}_{N} \hat{A} {\cal P}_{N}
$$
for an operator $\hat{A}$.
Then it is easy to see that (\ref{eqn:det1}) is written in the form
\begin{eqnarray}
\label{eqn:det1b}
\Bigg\{ \Psi_{N,T} ( {\bf f}; \vtheta ) 
\Bigg\}^2
&=& \det_{1 \leq i, j \leq N} 
\Bigg[ \delta_{ij}
+ \langle i | {\bf B} | \mu,x \rangle 
\langle \mu,x | {\bf C} |j \rangle \Bigg] \\
&=& \det_{1 \leq i, j \leq N}
\langle i | \Bigg[
\hat{\delta}_{N}+ {\bf B} {\bf C} \Bigg] | j \rangle,
\nonumber
\end{eqnarray}
where $\langle i | {\bf B} | \mu, x \rangle$
and $\langle \mu,x | {\bf C} | j \rangle$
are the two-dimensional row and
column vectors, respectively, given by
\begin{small}
\begin{eqnarray}
\langle i | {\bf B} | \mu, x \rangle 
&=&  \left( \begin{array}{lll}
\langle i | {\cal P}_{N}  
\hat{J}_{N}^{\bigtriangleup} {\cal P}_{N} 
|k \rangle \langle k | \hat{\chi}
\frac{1}{1-\hat{p}_{+} \hat{\chi}} | \mu, x \rangle  
 & &
-\langle i | {\cal P}_{N} 
\hat{J}_{N}^{\bigtriangleup}
{\cal P}_{N} 
|k \rangle \langle k |
\hat{J} \hat{p} \hat{\chi} | \mu, x \rangle \\
  & &
- \langle i |  {\cal P}_{N} \hat{J}_{N}^{\bigtriangleup}
{\cal P}_{N} | k \rangle \\
&& \quad \times \langle k | \hat{\chi}
\frac{1}{1-\hat{p}_{+} \hat{\chi}} | \nu, y \rangle 
\langle \nu, y | \hat{p} \hat{J} \hat{p}
\hat{\chi} | \mu, x \rangle 
\end{array}
\right), \nonumber\\
\langle \mu,x | {\bf C} | j \rangle
&=& \left( \matrix{
\langle \mu, x | \hat{p} \hat{J} {\cal P}_{N} | j \rangle \cr
\cr
-\langle \mu, x | \frac{1}{1-\hat{p}_{-} \hat{\chi}} {\cal P}_{N}
 | j \rangle } \right). \nonumber
\end{eqnarray}
\end{small}

The determinant (\ref{eqn:det1b}) is written using
the Fredholm determinant,
\begin{small}
\begin{eqnarray}
&& {\rm Det} \langle \mu,x | \Bigg[
I_{2}+{\bf C B} \Bigg] |\nu, y \rangle \nonumber\\
&=& {\rm Det} \left[ 
\begin{array}{ll}
\langle \mu,x | \nu, y \rangle 
 &
-\langle \mu,x | \hat{p} \hat{J}
|i \rangle \langle i| {\cal P}_{N} \hat{J}_{N}^{\bigtriangleup} 
{\cal P}_{N} 
|j \rangle \langle j | \hat{J} \hat{p} \hat{\chi} 
| \nu, y \rangle \\
+ \langle \mu,x | \hat{p} \hat{J} 
|i \rangle 
    & - \langle \mu,x | \hat{p} \hat{J} 
    |i \rangle \langle i |{\cal P}_{N} \hat{J}_{N}^{\bigtriangleup} 
    {\cal P}_{N} 
    |j \rangle \langle j |
    \hat{\chi}
    \frac{1}{1-\hat{p}_{+}\hat{\chi}}| \rho, z \rangle \\
  \  \times 
  \langle i | {\cal P}_{N} \hat{J}_{N}^{\bigtriangleup} {\cal P}_{N}
  |j \rangle \langle j |
\hat{\chi}
\frac{1}{1-\hat{p}_{+}\hat{\chi}}| \nu, y \rangle 
    & \qquad \qquad \times
   \langle \rho, z | \hat{p} \hat{J} \hat{p} \hat{\chi}
   | \nu, y \rangle \\
   &   \\
 -\langle \mu,x | \frac{1}{1-\hat{p}_{-}\hat{\chi}} 
 |i \rangle &
\langle \mu, x | \nu, y \rangle \\
 \ \times
 \langle i |{\cal P}_{N} \hat{J}_{N}^{\bigtriangleup} {\cal P}_{N} 
 |j \rangle \langle j | \hat{\chi}
\frac{1}{1-\hat{p}_{+} \hat{\chi}}| \nu, y \rangle
& 
+ \langle \mu, x | \frac{1}{1-\hat{p}_{-}\hat{\chi}} 
|i \rangle \langle i| {\cal P}_{N} \hat{J}_{N}^{\bigtriangleup}
{\cal P}_{N} |j \rangle \langle j | \hat{J}
\hat{p} \hat{\chi} | \nu, y \rangle \\
  &   + \langle \mu, x | 
  \frac{1}{1-\hat{p}_{-} \hat{\chi}} 
  |i \rangle \langle i| {\cal P}_{N} \hat{J}_{N}^{\bigtriangleup}
   {\cal P}_{N} 
  |j \rangle \langle j | \hat{\chi}
  \frac{1}{1-\hat{p}_{+}\hat{\chi}} | \rho, z \rangle \\
  & \qquad \qquad \times
  \langle \rho,z | \hat{p} \hat{J} \hat{p} \hat{\chi}
  | \nu, y \rangle 
\end{array} \right]
\nonumber\\
&=& {\rm Det} \left( \left[
\begin{array}{ll}
\langle \mu,x | \rho, z \rangle 
 &
-\langle \mu,x | \hat{p} \hat{J} |i \rangle \\
+ \langle \mu,x | \hat{p} \hat{J} | i \rangle
 & \quad \times
 \langle i | {\cal P}_{N} \hat{J}_{N}^{\bigtriangleup}
{\cal P}_{N} |j \rangle \langle j | \hat{J} \hat{p} \hat{\chi} 
| \rho, z \rangle 
\cr
\ \times
\langle i | {\cal P}_{N} \hat{J}_{N}^{\bigtriangleup} {\cal P}_{N} 
|j \rangle \langle j | \hat{\chi}
\frac{1}{1-\hat{p}_{+}\hat{\chi}}| \rho, z \rangle 
 &
+ \langle \mu, x | \hat{p} \hat{J} \hat{p} \hat{\chi}
| \rho, z \rangle \cr
   &   \cr
-\langle \mu,x | \frac{1}{1-\hat{p}_{-}\hat{\chi}} 
|i \rangle & \langle \mu, x | \rho, z \rangle
\cr
\ \times \langle i | {\cal P}_{N} \hat{J}_{N}^{\bigtriangleup} 
{\cal P}_{N} 
|j \rangle \langle j | \hat{\chi}
\frac{1}{1-\hat{p}_{+}\hat{\chi}}| \rho, z \rangle
&
+ \langle \mu, x | \frac{1}{1-\hat{p}_{-} \hat{\chi}} 
|i \rangle \cr
& \ \times
\langle i | {\cal P}_{N} \hat{J}_{N}^{\bigtriangleup}
{\cal P}_{N} |j \rangle \langle j | \hat{J} 
\hat{p} \hat{\chi} | \rho, z \rangle 
\end{array}
\right] \right. \nonumber\\
  && \hskip 8cm \times \left.
  \left[ \matrix{
  & \cr
  \langle \rho,z | \nu, y \rangle & 
  -\langle \rho,z |  \hat{p} \hat{J} \hat{p} \hat{\chi} 
  | \nu, y \rangle 
  \cr  & \cr & \cr
  0 & \langle \rho,z | \nu, y \rangle \cr
  & } \right] \right)
\nonumber\\
&=& 
{\rm Det} \left[
\begin{array}{ll}
\langle \mu,x | \nu, y \rangle 
 &
-\langle \mu,x | \hat{p} \hat{J} |i \rangle \\
+ \langle \mu,x | \hat{p} \hat{J} | i \rangle
 & \quad \times
 \langle i | {\cal P}_{N} \hat{J}_{N}^{\bigtriangleup}
{\cal P}_{N} |j \rangle \langle j | \hat{J} \hat{p} \hat{\chi} 
| \nu, y \rangle 
\cr
\ \times
\langle i | {\cal P}_{N} \hat{J}_{N}^{\bigtriangleup}
 {\cal P}_{N} 
 |j \rangle \langle j | \hat{\chi}
\frac{1}{1-\hat{p}_{+}\hat{\chi}}| \nu, y \rangle 
 &
+ \langle \mu, x | \hat{p} \hat{J} \hat{p} \hat{\chi}
| \nu, y \rangle \cr
   &   \cr
-\langle \mu,x | \frac{1}{1-\hat{p}_{-}\hat{\chi}} 
|i \rangle & \langle \mu, x | \nu, y \rangle
\cr
\ \times \langle i | {\cal P}_{N} \hat{J}_{N}^{\bigtriangleup}
{\cal P}_{N} 
|j \rangle \langle j | \hat{\chi}
\frac{1}{1-\hat{p}_{+}\hat{\chi}}| \nu, y \rangle
&
+ \langle \mu, x | \frac{1}{1-\hat{p}_{-} \hat{\chi}} 
|i \rangle \cr
& \ \times
\langle i | {\cal P}_{N} \hat{J}_{N}^{\bigtriangleup}
{\cal P}_{N} |j \rangle \langle j | \hat{J} 
\hat{p} \hat{\chi} | \nu, y \rangle 
\end{array}
\right], \nonumber
\end{eqnarray}
\end{small}
where $I_{2}$ denotes the unit matrix with size 2.
It is further rewritten as
\begin{small}
\begin{eqnarray}
&& {\rm Det} \langle \mu, x | \left( I_{2}
+ \left[ \matrix{ 1 & 0 \cr & \cr
0 & \frac{1}{1- \hat{p}_{-} \hat{\chi}} } \right] \right.
\nonumber\\
&& \qquad \qquad \qquad \times
\left[ \matrix{ \hat{p} \hat{J} 
|i \rangle \langle i | {\cal P}_{N} \hat{J}_{N}^{\bigtriangleup}
{\cal P}_{N}
|j \rangle \langle j |
\hat{\chi} &
\Big(-\hat{p} \hat{J} 
|i \rangle \langle i | {\cal P}_{N} \hat{J}_{N}^{\bigtriangleup} 
{\cal P}_{N} 
|j \rangle \langle j | \hat{J} 
   \hat{p}
  +\hat{p} \hat{J} \hat{p} \Big) \hat{\chi} \cr
    & \cr
-|i \rangle \langle i | {\cal P}_{N} \hat{J}_{N}^{\bigtriangleup}
 {\cal P}_{N} |j \rangle \langle j | \hat{\chi} 
& |i \rangle \langle i | {\cal P}_{N} \hat{J}_{N}^{\bigtriangleup}
{\cal P}_{N} 
|j \rangle \langle j | 
\hat{J} \hat{p} \hat{\chi} } \right] 
\nonumber\\
&& \hskip 11cm \times \left.
 \left[ \matrix{
    \frac{1}{1-\hat{p}_{+} \hat{\chi}} & 0 \cr
   & \cr
   0 & 1 } \right] \right) |\nu, y \rangle
\nonumber\\
&=& 
{\rm Det} \langle \mu, x | \left( 
\left[ \matrix{ 1 & 0 \cr & \cr
0 & 1- \hat{p}_{-} \hat{\chi} } \right]
\left[ \matrix{
    1-\hat{p}_{+} \hat{\chi} & 0 \cr
   & \cr
   0 & 1 } \right] \right. \nonumber\\
&& \hskip 2cm \left.
+ \left[ \matrix{ \hat{p} \hat{J} 
|i \rangle \langle i | {\cal P}_{N} \hat{J}_{N}^{\bigtriangleup} 
{\cal P}_{N}
|j \rangle \langle j |
\hat{\chi} &
\Big(-\hat{p} \hat{J} 
|i \rangle \langle i | {\cal P}_{N} \hat{J}_{N}^{\bigtriangleup} 
{\cal P}_{N} 
|j \rangle \langle j | \hat{J} 
   \hat{p}
  +\hat{p} \hat{J} \hat{p} \Big) \hat{\chi} \cr
    & \cr
-|i \rangle \langle i | {\cal P}_{N} \hat{J}_{N}^{\bigtriangleup}
 {\cal P}_{N} |j \rangle \langle j | \hat{\chi} 
& |i \rangle \langle i | {\cal P}_{N} \hat{J}_{N}^{\bigtriangleup}
{\cal P}_{N} |j \rangle \langle j |
\hat{J} \hat{p} \hat{\chi} } \right]
\right)
|\nu, y \rangle
\nonumber\\
&=& {\rm Det} \langle \mu, x | \left(
I_{2}
+ \left[
\matrix{ \hat{p} \hat{J}
|i \rangle \langle i | {\cal P}_{N} \hat{J}_{N}^{\bigtriangleup}
 {\cal P}_{N}
|j \rangle \langle j |
-\hat{p}_{+} &
-\hat{p} \hat{J}
|i \rangle \langle i | {\cal P}_{N} \hat{J}_{N}^{\bigtriangleup} 
{\cal P}_{N} 
|j \rangle \langle j | \hat{J} 
   \hat{p}
  +\hat{p} \hat{J} \hat{p} \cr
    & \cr
-|i \rangle \langle i | {\cal P}_{N} \hat{J}_{N}^{\bigtriangleup}
 {\cal P}_{N} 
 |j \rangle \langle j |
& |i \rangle \langle i | {\cal P}_{N} \hat{J}_{N}^{\bigtriangleup}
{\cal P}_{N}
|j \rangle \langle j | \hat{J} \hat{p} -\hat{p}_{-} }
\right] \hat{\chi} \right)
|\nu, y \rangle.
\nonumber
\end{eqnarray}
\end{small}
Here we have used the facts that
$$
{\rm Det} \langle \mu, x | 
\left[ \matrix{ 1 & 0 \cr & \cr
0 & 1- \hat{p}_{-} \hat{\chi} } \right]
| \nu, y \rangle
=1,
$$
and
$$
{\rm Det} \langle \mu, x | 
\left[ \matrix{
    1-\hat{p}_{+} \hat{\chi} & 0 \cr
   & \cr
   0 & 1 } \right] | \nu, y \rangle =1,
$$
which are consequences of definitions
(\ref{eqn:p+-}) of the operators
$\hat{p}_{+}$ and $\hat{p}_{-}$.
Then we arrive at
\begin{equation}
\label{eqn:Psi2}
\Bigg\{ \Psi_{N,T}( {\bf f}; \vtheta)
\Bigg\}^2
={\rm Det} \left(
I_{2} \delta_{\mu \nu} \delta(x-y) + 
\left[ \matrix{ \tS^{\mu,\nu}(x,y) & \tI^{\mu,\nu}(x,y) \cr
D^{\mu,\nu}(x,y) & \tS^{\nu,\mu}(y,x) } \right] 
\chi_{\nu}(y) \right),
\end{equation}
where
\begin{eqnarray}
\label{eqn:DSI1}
D^{\mu,\nu}(x,y) 
&=& -\langle \mu,x | i \rangle 
\langle i | {\cal P}_{N} 
({\cal P}_{N} \hat{J}{\cal P}_{N})^{\bigtriangleup} {\cal P}_{N}
 | j \rangle
\langle j | \nu, y \rangle, \\
S^{\mu,\nu}(x,y) 
&=& \langle \mu, x | \hat{p} \hat{J}
| i \rangle \langle i | 
{\cal P}_{N}  
({\cal P}_{N} \hat{J} {\cal P}_{N})^{\bigtriangleup}
{\cal P}_{N} | j \rangle
\langle j | \nu, y \rangle,
\nonumber\\
I^{\mu,\nu}(x,y)
&=& -\langle \mu, x | \hat{p} \hat{J} | i \rangle
\langle i | {\cal P}_{N}  
({\cal P}_{N} \hat{J}{\cal P}_{N})^{\bigtriangleup}
 {\cal P}_{N} 
| j \rangle
\langle j | \hat{J} \hat{p} | \nu, y \rangle,
\nonumber
\end{eqnarray}
and
\begin{eqnarray}
\label{eqn:DSI2}
\\
\tS^{\mu,\nu}(x,y) &=& S^{\mu,\nu}(x,y)
-\langle \mu,x | \hat{p}_{+} | \nu, y \rangle \nonumber\\
&=& \left\{
\begin{array}{ll}
\langle \mu, x | \hat{p} \hat{J} |i \rangle
\langle i |{\cal P}_{N} 
({\cal P}_{N} \hat{J} {\cal P}_{N})^{\bigtriangleup}
{\cal P}_{N} | j \rangle
\langle j |\nu, y \rangle,
& \mbox{if $\mu \geq \nu$,} \\
& \\
-\langle \mu, x | \hat{p} \hat{J} | i \rangle
\langle i | 
\Bigg( \hat{J}^{\bigtriangleup}
- {\cal P}_{N}
({\cal P}_{N} \hat{J} {\cal P}_{N})^{\bigtriangleup}
{\cal P}_{N}\Bigg) | j \rangle \langle j | \nu, y \rangle,
& \mbox{if $\mu < \nu$,}
\end{array}
\right. \nonumber\\
\tI^{\mu,\nu}(x,y) &=& I^{\mu,\nu}(x,y)
+ \langle \mu, x | \hat{p} \hat{J} \hat{p} | \nu, y 
\rangle \nonumber\\
&=& \langle \mu, x | \hat{p} \hat{J} | i \rangle
\langle i | \Bigg( \hat{J}^{\bigtriangleup}
-{\cal P}_{N}
({\cal P}_{N}\hat{J}{\cal P}_{N})^{\bigtriangleup}
{\cal P}_{N} \Bigg) | j \rangle
\langle j | \hat{J} \hat{p} | \nu, y \rangle, \nonumber
\end{eqnarray}
where we have used the equality (\ref{eqn:equalityA})
for $\hat{A}=\hat{J}$.

In \cite{Rains00} Rains introduced Fredholm pfaffian,
denoted here by ${\rm PF}$, and proved a useful equality
\begin{equation}
\label{eqn:Rains1}
\Bigg\{ {\rm PF}(J_{2}+K) \Bigg\}^2
={\rm Det}(I_{2}+J_{2}^{-1} K),
\end{equation}
for any antisymmetric $2 \times 2$ matrix kernel $K$,
where 
\begin{equation}
\label{eqn:J2}
J_{2}=\left( \matrix{ 0 & 1 \cr -1 & 0} \right).
\end{equation}
Then (\ref{eqn:Psi2}) implies the Fredholm pfaffian
representation of the multi-time characteristic function
\begin{equation}
\label{eqn:Psi3}
\Psi_{N,T}( {\bf f}; \vtheta)
={\rm PF} \left(
J_{2} \delta_{\mu \nu} \delta(x-y) + 
\left[ \matrix{ D^{\mu,\nu}(x,y) & \tS^{\nu,\mu}(y,x) \cr
-\tS^{\mu,\nu}(x,y) & -\tI^{\mu,\nu}(x,y) }\right] 
\chi_{\nu}(y) \right).
\end{equation}

\setcounter{equation}{0}
\section{Pfaffian Process}

Let
\begin{equation}
\label{eqn:matrixA}
A^{\mu,\nu}(x,y)=\left[ \matrix{ D^{\mu,\nu}(x,y) & 
\tS^{\nu,\mu}(y,x) \cr
-\tS^{\mu,\nu}(x,y) & -\tI^{\mu,\nu}(x,y) } \right],
\end{equation}
and construct
$2 \sum_{\mu=1}^{M+1} N_{\mu} \times 2 \sum_{\mu=1}^{M+1} N_{\mu}$
antisymmetric matrices
$$
A \Bigg( \x_{N_{1}}^{(1)}, \x_{N_{2}}^{(2)}, \cdots,
\x_{N_{M+1}}^{(M+1)} \Bigg)
= \left( A^{\mu,\nu} \Bigg(x_{i}^{(\mu)}, x_{j}^{(\nu)} \Bigg)
\right)_{1 \leq i \leq N_{\mu}, 1 \leq j \leq N_{\nu}, 
1 \leq \mu, \nu \leq M+1},
$$
for $N_{m}=1,2, \cdots, N, 1 \leq m \leq M+1$.
By the definition of Fredholm pfaffian \cite{Rains00}
and the equality (\ref{eqn:Psi3}), we can establish the
following statement.
\begin{thm}
\label{thm:Main1}
The non-colliding system of Brownian motions $\X(t)$ is a 
pfaffian process in the sense that
any multi-time correlation function is given by a pfaffian
\begin{eqnarray}
\label{eqn:rho2}
&& \rho_{N,T} \left(t_{1}, \{\x^{(1)}_{N_1}\}; t_2, \{\x^{(2)}_{N_2}\}; 
\dots; t_{M+1}, \{\x^{(M+1)}_{N_{M+1}}\} \right) \\
&& \qquad ={\rm Pf} \left[
A \Bigg( \x_{N_{1}}^{(1)}, \x_{N_{2}}^{(2)}, \cdots,
\x_{N_{M+1}}^{(M+1)} \Bigg) \right].
\end{eqnarray}
\end{thm}

\noindent{\bf Remark.} \quad
The pfaffian processes considered here may be regarded as 
the continuous space-time version of the pfaffian
point processes and pfaffian Schur processes
studied by Sasamoto and Imamura \cite{SI04}
and by Borodin and Rains \cite{BR04}.
The processes studied in \cite{FNH99} are
also pfaffian processes, since the `quaternion
determinantal expressions', in the sense
of Dyson and Mehta \cite{Dyson70,Mehta89,Mehta91},
of correlation functions
are readily transformed to pfaffian expressions.

Let $H_{i}(x)$ be the $i$-th Hermite polynomial
\begin{eqnarray}
  H_{i}(x) &=& e^{x^2} \left(- \frac{d}{dx} \right)^{i} e^{-x^2}
  \nonumber\\
    &=& i ! \sum_{j=0}^{[i/2]} (-1)^{j}
    \frac{(2x)^{i-2j}}{j ! (i-2j)!},
\nonumber
\end{eqnarray}
where $[a]$ denotes the greatest integer not greater than $a$.
The Hermite polynomials satisfy the orthogonal relations
\begin{equation}
\label{eqn:orth}
 \int_{\R} dx \, e^{-x^2} H_{i}(x) H_{j}(x)= 2^{i} i !
 \sqrt{\pi} \delta_{ij}, \qquad
 i, j =0,1,2, \cdots.
\end{equation}
Set
$$
  c_{1}=\sqrt{\frac{t_{1}(2T-t_{1})}{T}}, \qquad
  z_{1}=\sqrt{\frac{2T-t_{1}}{t_{1}}},
$$
and
\begin{equation}
\alpha_{i j}=\left\{
\begin{array}{ll}
2^{-i}c_1^{i}\delta_{i j},
& \mbox{if $i$ is even,}
\\ 
2^{-i}c_1^{i}
\Big\{ \delta_{i j}- 2(i-1)\delta_{i-2 \, j} \Big\},
& \mbox{if $i$ is odd.}
\end{array}
\right.
\label{def:alpha}
\end{equation}

Now we specify polynomials $\{M_{i}(x)\}$ as
\begin{equation}
\label{eqn:M}
  M_{i}(x)= b_{i} z_{1}^{-i} \sum_{j=0}^{i} \alpha_{i j}
  H_{j} \left( \frac{x}{c_{1}} \right) z_{1}^{j},
  \quad i=0,1,2, \cdots 
\end{equation}
with $b_{i}=\Big\{r_{[i/2]} \Big\}^{-1/2}$, where \,
$\displaystyle{
  r_{i}=\frac{1}{\pi} \Gamma(i+1/2) \Gamma(i+1)
  \left(\frac{t_{1}^2}{T} \right)^{2i+1/2}.
}$
Set
$$
    J_{N}=I_{N/2} \otimes J_{2},
$$
where 
$I_{N/2}$ denotes the unit matrix
with size $N/2$ and $J_{2}$ is given by (\ref{eqn:J2}),
and let $J$ be the semi-infinite matrix obtained as the
$N \to \infty$ limit of $J_{N}$.
By the orthogonality of Hermite polynomials
(\ref{eqn:orth}), 
we can show through (\ref{eqn:R}) 
with the choice (\ref{eqn:M}) that \cite{NKT03,KNT04}
$$
   \langle i | \hat{J} | j \rangle
= \langle i | \mu, x \rangle \langle \mu, x |
\hat{J} | \nu, y \rangle \langle \nu, y | j \rangle
= J_{ij}, \qquad
i, j =1,2, \cdots.
$$
Since $J_{N}^2=-I_{N}$ for any even $N \geq 2$,
this implies that the matrix 
$\Big( \langle i | \hat{J} | j \rangle \Big)_{i,j=1,2,\cdots}$
is invertible and 
\begin{eqnarray}
\label{eqn:A-1}
 \langle i | \hat{J}^{\bigtriangleup} | j \rangle &=& 
  = - J_{ij}, \\
\label{eqn:A-2}
\langle i | {\cal P}_{N}
({\cal P}_{N} \hat{J}{\cal P}_{N})^{\bigtriangleup}{\cal P}_{N}
|j \rangle 
&=& \langle i | {\cal P}_{N} \hat{J}^{\bigtriangleup}
{\cal P}_{N} 
| j \rangle \\
&=& 
\left\{
\begin{array}{ll}
-J_{ij}=-(J_{N})_{ij},
& \mbox{if $1 \leq i, j \leq N $,}
\\ 
0,
& \mbox{otherwise}
\end{array}
\right. \nonumber
\end{eqnarray}
for $i,j =1,2, \cdots.$

If we write
\begin{eqnarray}
&& \langle \mu,x | i \rangle 
= b_{i} R_{i-1}^{(\mu)}(x), \nonumber\\
&& \langle i | \hat{J} \hat{p} | \mu, x \rangle
= -\langle \mu, x | \hat{p} \hat{J} | i \rangle 
= b_{i} \Phi_{i-1}^{(\mu)}(x), \nonumber
\end{eqnarray}
$i=1,2, \cdots, \mu=1,2, \cdots, M+1, x \in \R$, then
the functions (\ref{eqn:DSI1}) are written as
\begin{eqnarray}
\label{eqn:DSI3}
D^{\mu,\nu}(x,y) &=& \sum_{i=0}^{N/2-1} \frac{1}{r_{i}}
\Bigg[ R_{2i}^{(\mu)}(x) R_{2i+1}^{(\nu)}(y)
   - R_{2i+1}^{(\mu)}(x) R_{2i}^{(\nu)}(y) \Bigg], \\
S^{\mu,\nu}(x,y) &=& \sum_{i=0}^{N/2-1} \frac{1}{r_{i}}
\Bigg[ \Phi_{2i}^{(\mu)}(x) R_{2i+1}^{(\nu)}(y)
   - \Phi_{2i+1}^{(\mu)}(x) R_{2i}^{(\nu)}(y) \Bigg], \nonumber\\
I^{\mu,\nu}(x,y) &=& - \sum_{i=0}^{N/2-1} \frac{1}{r_{i}}
\Bigg[ \Phi_{2i}^{(\mu)}(x) \Phi_{2i+1}^{(\nu)}(y)
   - \Phi_{2i+1}^{(\mu)}(x) \Phi_{2i}^{(\nu)}(y) \Bigg], \nonumber
\end{eqnarray}
and Equations (\ref{eqn:DSI2}) become
\begin{eqnarray}
\label{eqn:DSI4}
\qquad
\tS^{\mu, \nu}(x,y) 
&=& \left\{
\begin{array}{ll}
\displaystyle{\sum_{i=0}^{N/2-1} \frac{1}{r_{i}}
\Bigg[ \Phi_{2i}^{(\mu)}(x) R_{2i+1}^{(\nu)}(y)
   - \Phi_{2i+1}^{(\mu)}(x) R_{2i}^{(\nu)}(y) \Bigg]},
& \mbox{if $\mu \geq \nu$,} \\ 
& \\
\displaystyle{-\sum_{i=N/2}^{\infty} \frac{1}{r_{i}}
\Bigg[ \Phi_{2i}^{(\mu)}(x) R_{2i+1}^{(\nu)}(y)
   - \Phi_{2i+1}^{(\mu)}(x) R_{2i}^{(\nu)}(y) \Bigg]},
& \mbox{if $\mu < \nu$,}
\end{array}
\right. \\
\tI^{\mu,\nu}(x,y) &=&
\sum_{i=N/2}^{\infty} \frac{1}{r_{i}}
\Bigg[ \Phi_{2i}^{(\mu)}(x) \Phi_{2i+1}^{(\nu)}(y)
   - \Phi_{2i+1}^{(\mu)}(x) \Phi_{2i}^{(\nu)}(y) \Bigg]. \nonumber
\end{eqnarray}

Theorem \ref{thm:Main1} with the expressions (\ref{eqn:DSI3}),
(\ref{eqn:DSI4}) of the
elements of matrix (\ref{eqn:matrixA}) is
equivalent with Theorem 3 reported in
\cite{KNT04}, although the latter was given in the form
of quaternion determinant.
The present argument will be generalized to discuss
other non-colliding systems of diffusion particles
reported in \cite{KT04a,KT05}.

\section*{Acknowledgment}
This manuscript is based on the joint work with
H. Tanemura.
The author would like to thank T. Sasamoto 
for useful discussions on the papers
\cite{Rains00} and \cite{BR04}.


\end{document}